# Joint optimal power allocation and sensing threshold selection for SU's capacity maximization in SS CRNs

F. T. Foukalas, P. T. Mathiopoulos and G. T. Karetsos

We propose a joint optimal power allocation and sensing threshold selection for capacity maximization at the secondary user (SU) in spectrum sharing (SS) cognitive radio networks (CRNs). Hence, both optimal power allocation (OPA) and spectrum sensing (SpSe) is considered in the SS CRNs model. The obtained results show that such a joint optimal selection improves the performance of the SU by maximizing its capacity.

*Introduction:* In spectrum sharing (SS) cognitive radio networks (CRNs), optimal power allocation (OPA) and spectrum sensing (SpSe) are used for the protection of the primary user (PU) from harmful interference caused by the secondary user (SU). Furthermore, for capacity maximization of the SU, the main parameters related to OPA, i.e. the SU's transmit power, $P_t$, is adapted according to the received signal-to-noise ratio (SNR), $\gamma_s$, and related to SpSe, i.e. sensing threshold, $\eta$, and sensing time, $\tau$, for a given sensed SNR, $\gamma$, need to be also carefully selected [1]. Previous studies on SU's capacity maximization include SS CRNs models with SpSe [5],[6] or without SpSe [2],[3]. For the former and more general case, the optimization presented in [5] is considered over $P_t$ and $\tau$ assuming $\eta$ to be constant. In [6] although the effects of $\eta$ as a variable are studied, the research is focused on the interference caused to the PU rather than the optimization of the SU's capacity. Thus, a more general approach is presented in this letter where a jointly OPA and SpSe threshold selection is considered so that the SU's capacity is maximized over $P_t$ and $\eta$.

*SU's capacity maximization:* Following [5], we consider a SS CRN system model with OPA and SpSe capabilities operating in a combination of an additive white Gaussian noise (AWGN) channel with mean zero and variance $N_0$ and Rayleigh fading [4]. Its SU capacity can be obtained as [5]

$$C_s = C_0 \pi_0 (1 - p_f) + C_1 \pi_0 p_f + C_1 \pi_1 p_d + C_0 \pi_1 (1 - p_d) \tag{1}$$

where $\pi_0$ and $\pi_1$ are the probabilities for which the PU is idle or active, respectively, $p_d$ denotes the probability of correct detection and $p_f$ the probability of false alarm for SpSe. $C_0$ is the SU's capacity assuming that the PU correctly or falsely is not detected and $C_1$ is the SU's capacity when the PU is correctly or falsely detected [4], in which cases the SU's OPA are $P_t^0$ and $P_t^1$, respectively. Maximization of $C_s$ over $P_t^0, P_t^1$ and $\eta$ will be performed assuming $P_t^0 > P_t^1$ for the protection of PU under two constraints. The first relates to the SU average power, $P_{av}$, as follows:

$$\overline{P} \leq P_{av} \tag{C1}$$

where $\overline{P} = \pi_0 E(P_t^0)(1 - p_f) + \pi_0 E(P_t^1) p_f + \pi_1 E(P_t^1) p_d + \pi_1 E(P_t^0)(1 - p_d)$ is the average transmit power of SU with $E(\cdot)$ being the expectation over the probability density function (pdf) of the Rayleigh fading channel. The second relates to the peak interference power, $I_{pk}$, at the PU as follows

$$G_{sp}\left[\pi_0 P_t^0 (1 - p_f) + \pi_0 P_t^1 p_f + \pi_1 P_t^1 p_d + \pi_1 P_t^0 (1 - p_d)\right] \leq I_{pk} \tag{C2}$$

where $G_{sp}$ is the channel power gain of the link between the SU and PU. Furthermore, since for SpSe the optimal sensing threshold, $\eta^*$, is related to the optimal probability of detection, $p_d^*$, at which the PU transmission is protected from the SU's transmission, a third constraint must be used

$$p_d(\eta, \tau) \geq p_d^* \qquad (C3)$$

where $\tau$ is assumed to be constant.

Thus the problem of $C_s$ maximization can be stated as

$$\max_{\{P_t^0, P_t^1, \eta\}} C_s \qquad (2)$$

subject to (C1), (C2), and (C3)

Solving this joint optimization problem will lead to the optimal pair of values for the cut-off value on the received SNR at SU, $\gamma_s^*$, for the OPA and the optimum sensing threshold value, $\eta^*$, for the SpSe.

*Joint optimal power allocation and sensing threshold selection*: Since maximization of $C_s$ in (2) is a joint OPA and SpSe optimization problem, an iterative algorithm should be used. In particular, we rely on the sub-gradient method since: *a)* the objective function is complex non-differentiable, and *b)* the problem is convex over power and over sensing threshold since for $\eta > 0$, $p_f' > 0$ and $p_d' > 0$ is applied [3],[4]. Thus, it can be concluded that both $p_f$ and $p_d$ are concave up on $\eta$ [7]. The sub-gradient can be expressed from (C1) as

$$g = P_{av} - \overline{P} \qquad (3)$$

Using (3), we obtain the cut-off value $\gamma_s^*$ for the OPA. This value specifies the SU's OPA policies for a $P_{av}$ denoted in [3] as $P_t^0 / P_{av}$ and $P_t^1 / P_{av}$ obtained without and with interference power constraint, respectively. The updates of $\gamma_s^*$ are repeated until the convergence rule is reached [8]. The iterative updates are performed over the transmit powers $P_t^0, P_t^1$, that results in the cut-off value in SNR $\gamma_s^*$, for a specific $\eta$ where the constraint, $p_d^*$, determines the optimal sensing threshold, $\eta^*$. Thus, the optimum pair of values $[\gamma_s^*, \eta^*]$ can be obtained.

*Numerical results*: Figure 1 illustrates the performance of $C_s$ obtained from the joint optimization problem in (2), versus $\eta$, for different values of $\gamma$ and $P_{av}$. As in [3] and [5], for the performance evaluation results we have assumed that for Rayleigh fading channels the channel power gains (exponentially distributed) are assumed with unit mean, AWGN with variance $N_0 = 1$ and $\pi_1 = 0.4$. Furthermore, for the OPA, the constraint on peak interference power is assumed to be $I_{pk} = 0dB$ while for the SpSe, $\tau = 1ms$ [4]. The performance evaluation results obtained clearly show that $C_s$ increases as $\gamma$ decreases and/or $P_{av}$ increases while its improvement becomes negligible when $P_{av} < I_{pk}$.

Figure 2 illustrates the throughput computed as $\xi_s = (T - \tau / T) C_s$ based on the SpSe and frame transmission models that have been proposed in [5], where $T$ is the frame duration, $T - \tau$ is the frame duration for data transmission and $C_s$ is taken from (1). Therefore $\xi_s$ represents the transmitted bits per frame, is the performance metric at the secondary link and it is maximized over the sensing time $\tau$ for different optimal probabilities of detection $p_d^*$ using the maximization in (2). The performance results have been obtained for $T = 100ms$, $P_{av} = 15dB$, $I_{pk} = 0dB$ and $\gamma = -10db$. Furthermore, different target values of $p_d^*$ are assumed that correspond to specific SpSe thresholds, $\eta^*$. However, these $\eta^*$ values are identical for each target value $p_d^*$ as depicted in Fig. 2. This is reasonable since the maximization problem is assumed over $\tau$ and not over $\eta$. This also shows that a proper selection of $\eta$ for the SpSe and $P_t$ for the OPA provides an additional $C_s$ maximization to the one achieved by the joint optimization over $P_t$ and $\tau$.

*Conclusion*: A joint OPA and SpSe threshold selection for capacity maximization at the SU in SS CRNs has been proposed. This joint optimization leads to further capacity maximization as compared to the one achieved by joint optimization over the transmit power and sensing time. This maximization has also shown that the capacity can be further improved by properly selecting the SpSe threshold based on the sensed SNR.

**Authors' affiliations:**
F.Foukalas* (Dept. of Informatics and Telecommunications, National Kapodistrian University of Athens, Ilisia, Athens, Greece), P. T. Mathiopoulos (Institute for Space Applications and Remote Sensing, National Observatory of Athens, Metaxa, Athens, Greece) and G.Karetsos (Dept. of Information Technology and Telecommunications, TEI of Larissa, Larissa, Greece)

* Corresponding author (foukalas@di.uoa.gr)

**Figure captions:**

**Fig. 1** Capacity $C_s$ vs. sensing threshold $\eta$ for different $\gamma$ and $P_{av}$ assuming $I_{pk} = 0dB$

**Fig. 2** Throughput $\xi_s$ vs. sensing time, $\tau$ for different optimal $p_d^*$ with $P_{av} = 15dB$ and $I_{pk} = 0dB$

Figure 1

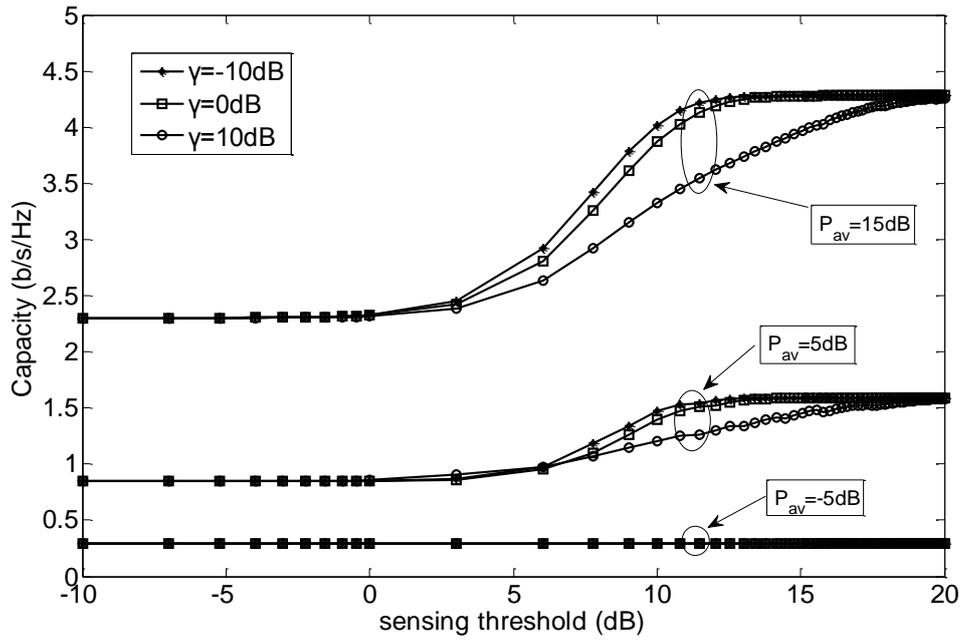

Figure 2

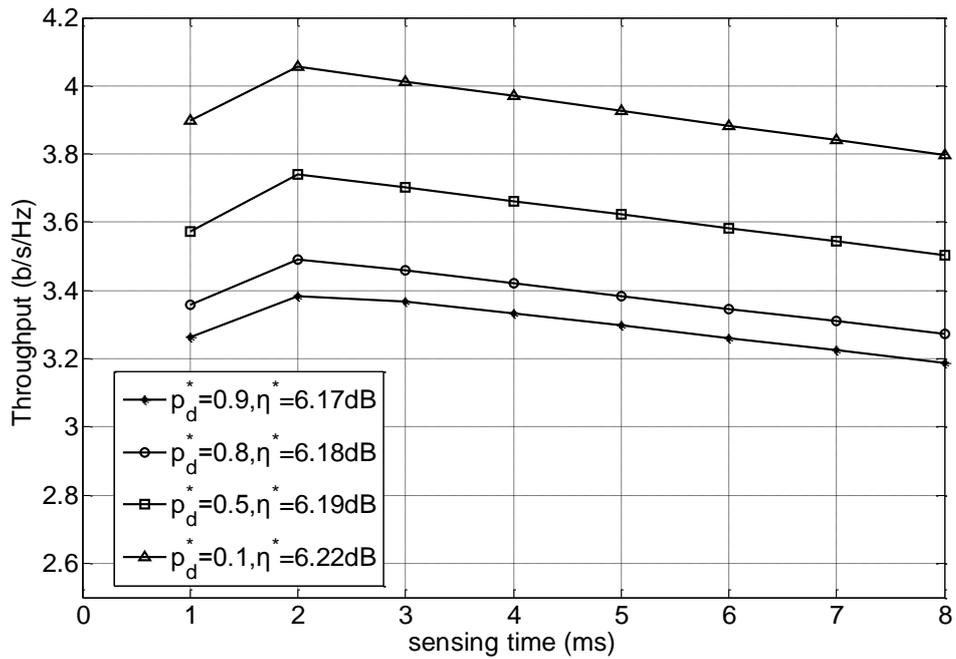